\documentclass[10pt,conference]{IEEEtran}
\usepackage{color}

\newtheorem{thm}{Theorem}

\newtheorem{prop}[thm]{Proposition}

\newtheorem{defn}[thm]{Definition}

\ifCLASSINFOpdf
   \usepackage[pdftex]{graphicx}
\else
\fi
\usepackage{amssymb}
\begin{document}
%
\title{Recovery guarantees for compressed sensing with unknown errors}

\author{\IEEEauthorblockN{Simone Brugiapaglia\IEEEauthorrefmark{1}, 
Ben Adcock\IEEEauthorrefmark{2}}
\IEEEauthorblockA{Department of Mathematics\\
Simon Fraser University\\
Burnaby, BC V5A 1S6, Canada\\
Email: 
\IEEEauthorrefmark{1}simone\_brugiapaglia@sfu.ca, 
\IEEEauthorrefmark{2}ben\_adcock@sfu.ca}
\and 
\IEEEauthorblockN{Richard K. Archibald}
\IEEEauthorblockA{Computer Science and Mathematics Division\\
Oak Ridge National Laboratory\\
Oak Ridge, TN 37831, USA\\
Email: archibaldrk@ornl.gov}
}


%


\maketitle
\thispagestyle{plain}
\pagestyle{plain}

\begin{abstract}
From a numerical analysis perspective, assessing the robustness of $\ell^1$-minimization is a fundamental issue in compressed sensing and sparse regularization. Yet, the recovery guarantees available in the literature usually depend on \emph{a priori} estimates of the noise, which can be very hard to obtain in practice, especially when the noise term also includes unknown discrepancies between the finite model and data. In this work, we study the performance of $\ell^1$-minimization when these estimates are not available, providing robust recovery guarantees for quadratically constrained basis pursuit and random sampling in bounded orthonormal systems. Several applications of this work are approximation of high-dimensional functions, infinite-dimensional sparse regularization for inverse problems, and fast algorithms for non-Cartesian Magnetic Resonance Imaging.
\end{abstract}

\IEEEpeerreviewmaketitle

\section{Introduction}

In \emph{Compressed Sensing} (CS) and sparse representations  we deal with underdetermined linear systems of equations 
\begin{equation}
\label{eq:noisymeas}
\mathbf{y} = \mathbf{A} \mathbf{x} + \mathbf{n},
\end{equation}
where $\mathbf{A}\in\mathbb{C}^{m\times N}$, with $m \ll N$, is the sensing matrix,  $\mathbf{x}\in\mathbb{C}^N$ is an unknown signal, and $\mathbf{y}\in\mathbb{C}^m$ is the vector of measurements perturbed by noise $\mathbf{n}\in\mathbb{C}^m$ \cite{Candes2006,Donoho2006}. This corruption could be due to physical noise produced by the measuring device, to approximation errors in the model, or to numerical factors. Some examples are model error in inverse problems such as MRI \cite{Guerquin-Kern2012,Adcock2016C}, the expansion error in infinite-dimensional CS when truncating the signal to its finite dimensional representation \cite{Adcock2016A,Adcock2016}, or the quadrature error involved in the evaluation of the bilinear form associated with a PDE \cite{Brugiapaglia2015,Brugiapaglia2016,PhDThesisSimone}.

\begin{figure}
\centering
\includegraphics[width = 7cm]{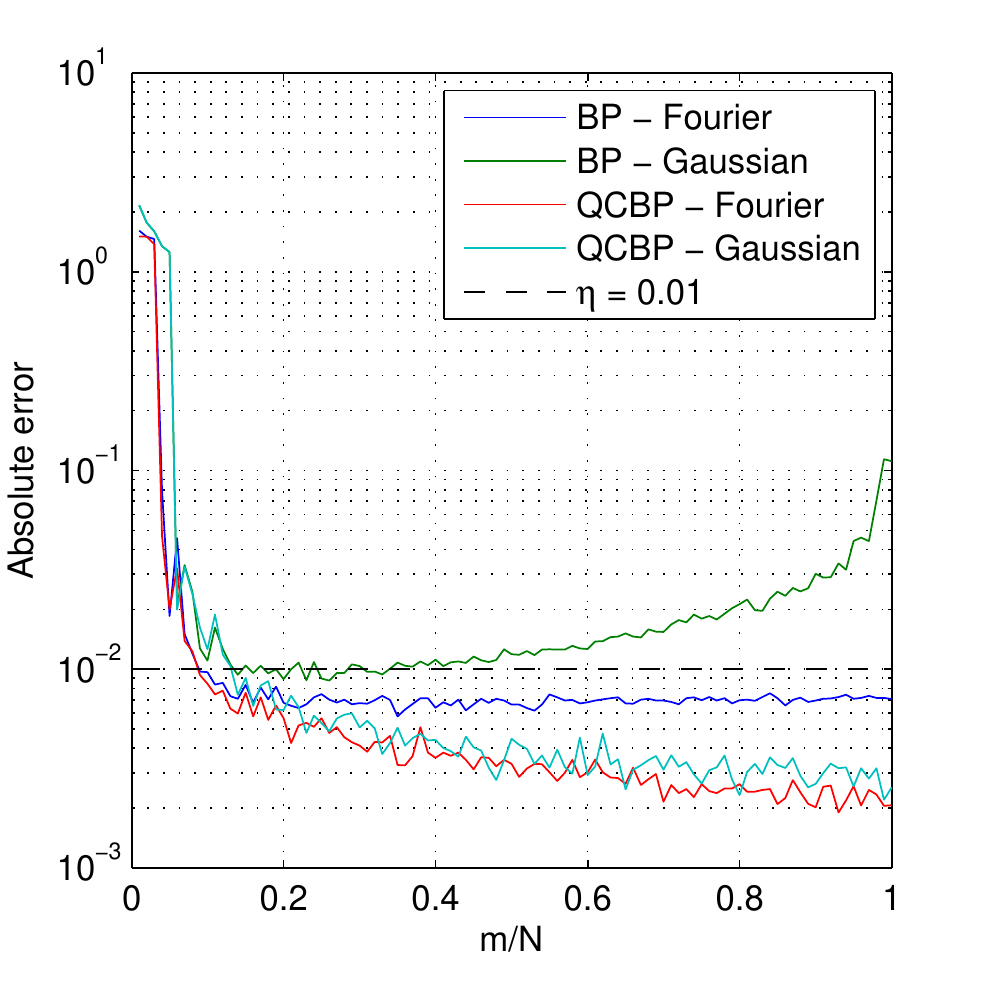}
\caption{\label{fig:FouriervsGauss}Numerical assessment of  BP ($\eta = 0$) and QCBP (with $\eta = 0.01$) for Fourier and Gaussian measurements corrupted by noise $\mathbf{n}$ of magnitude $\|\mathbf{n}\|_2=0.01$. The solution $\mathbf{x}$ is a randomly generated $10$-sparse vector in $\mathbb{C}^{1000}$. The absolute error $\|\mathbf{x}-\hat{\mathbf{x}}(\eta)\|_2$ is plotted as a function of the ratio $m/N$. The results are produced using the MATLAB package SPGL1 \cite{vanderBerg2008}. The QCBP solver, where relation (\ref{eq:noisebound}) holds, is very robust for both Fourier and Gaussian measurements. On the contrary, for BP, where relation (\ref{eq:noisebound}) does not hold anymore, the situation is different: the solver's performance highly depends on the type of measurements and on the ratio $m/N$.  Notably, Fourier measurements, coming from randomly subsampling the rows of the DFT matrix, are much more stable in the BP case when $m/N \rightarrow 1$ than Gaussian measurements. }
\end{figure}

A standard tool to regularize the inverse problem (\ref{eq:noisymeas}) and recover a good approximation $\hat{\mathbf{x}}(\eta)$ to the solution $\mathbf{x}$ (assumed to be sparse or compressible) is the \emph{Quadratically Constrained Basis Pursuit} (QCBP) optimization program
\begin{equation}
\label{eq:QCBP}
\hat{\mathbf{x}}(\eta)\in\arg\min_{\mathbf{z}\in\mathbb{C}^N} \|\mathbf{z}\|_1, \quad s.t. \quad \|\mathbf{A} \mathbf{z} - \mathbf{y}\|_2 \leq \eta,
\end{equation}
also called \emph{Basis Pursuit} (BP) when $\eta = 0$.
Usually, in order to study the recovery guarantees of (\ref{eq:QCBP}), the parameter $\eta$ is assumed to control the noise   magnitude, i.e.,
\begin{equation}
\label{eq:noisebound}
\|\mathbf{n}\|_2 \leq \eta.
\end{equation}
Indeed, under the regime (\ref{eq:noisebound}) and with suitable hypotheses on the sensing matrix $\mathbf{A}$ (e.g., based on the \emph{restricted isometry property}), the following type of recovery error estimate holds with high probability  
\begin{equation}
\label{eq:QCBPerror}
\|\mathbf{x} - \hat\mathbf{x}(\eta)\|_2 \lesssim \frac{\sigma_s(\mathbf{x})_1}{\sqrt{s}}  + \eta, 
\end{equation}
where $\sigma_s(\mathbf{x})_1$ is the best $s$-term approximation error of $\mathbf{x}$ with respect to the $\ell^1$-norm  \cite{Cohen2009,Foucart2013}. Unfortunately, \emph{a priori} estimates of the noise of the form (\ref{eq:noisebound}) may not be available in real applications of CS. Moreover, since the recovery error estimate (\ref{eq:QCBPerror}) is sensitive to $\eta$, the choice of this parameter is crucial  (see Figure~\ref{fig:FouriervsGauss}). In practice, one could resort to cross-validation in order to tune this parameter, but this technique could be time-consuming or inaccurate and it is not properly understood from a theoretical perspective \cite{Doostan2011}. 

The goal of this work is to establish robust recovery guarantees for QCBP (and BP)  under the regime 
\begin{equation}
\label{eq:noisegtreta}
\|\mathbf{n}\|_2 \geq \eta.
\end{equation}
In this scenario, recovery estimates analogous to (\ref{eq:QCBPerror})--where $\eta$ is replaced by $\|\mathbf{n}\|_2$--hold for BP \cite{Foucart2013}. They are based on the so-called quotient property, which is known to be satisfied only by random Gaussian matrices \cite{Wojtaszczyk2010} and by Weibull matrices \cite{Foucart2014}, under suitable restrictions on the number of measurements $m$. Similar robust recovery estimates are also available for algorithms such as iterative hard thresholding, CoSaMP, and orthogonal matching pursuit \cite{Foucart2013}. Yet, all these techniques require an \emph{a priori} knowledge of the sparsity level $s$ that is not necessary for QCBP.

Here, we prove robust recovery error estimates for QCBP (and BP) when the matrix $\mathbf{A}$ is built by random sampling from \emph{bounded orthonormal systems} \cite{Rauhut2010}. In particular,  under suitable hypotheses involving the restricted isometry constants and the singular values of $\mathbf{A}$, we provide  recovery error estimates in probability of the form (see Theorem~\ref{thm:robustrecovery})
\begin{equation}
\label{eq:QCBProbust}
\|\mathbf{x} - \hat\mathbf{x}(\eta)\|_2 
\lesssim \frac{\sigma_s(\mathbf{x})_1}{\sqrt{s}}  
+  \eta+ \mathcal{L}^{\frac12} \max \left \{ \|\mathbf{n}\|_2 - \eta , 0 \right \}, 
\end{equation}
where the factor $\mathcal{L}$ is polylogarithmic in $N$ and $s$, and can be defined, for example, as in (\ref{eq:defL}). The effect of the unknown error $\mathbf{n}$ is encapsulated in the third term on the right-hand side (compare with (\ref{eq:QCBPerror})).  As is to be expected, this term approaches zero as the estimation of the model error $\eta$ improves.

This analysis has several potential extensions and applications.  First, the case of weighted $\ell^1$-minimization.  This is used notably in high-dimensional function approximation and interpolation \cite{Adcock2016B,RauhutWard,WebsterPoly}, with applications in uncertainty quantification for parametric PDEs.  Second, fast methods for non-Cartesian MRI, where model error arises from gridding non-uniform Fourier data to a uniform grid, and can seriously hamper reconstruction quality \cite{Adcock2016D,Guerquin-Kern2012}.  Third, low-rank matrix recovery.

\section{Tools from compressed sensing}
We first recall some concepts from CS that constitute the foundations that our analysis will be built upon: the null space and the restricted isometry properties (Section~\ref{sec:NSP}), and random sampling in bounded orthonormal systems (Section~\ref{sec:BOS}). 

In the following, for every $k \in \mathbb{N}$, we define $[k]:=\{1,\ldots,k\}$ and $[k]_0:=\{0,\ldots,k-1\}$. Moreover, we denote the set of $s$-sparse vectors in $\mathbb{C}^N$ as $\Sigma_s^N$. 

\subsection{Robust null space and restricted isometry properties}
\label{sec:NSP}
The first tool involved in our analysis is the so-called \emph{$\ell^q$-robust Null Space Property} (NSP) \cite[Chapter 4]{Foucart2013}. 
\begin{defn}[$\ell^q$-robust null space property]
\label{def:NSP}
Given $q \geq 1$, the matrix $\mathbf{A} \in\mathbb{C}^{m \times N}$ satisfies the \emph{$\ell^q$-robust null space property} of order $s$ (with respect to the norm $\|\cdot\|$) with constants $0 < \rho < 1$ and $\tau >0$ if, for any set $S\subseteq [N]$ with $|S|\leq s$, it holds
\begin{equation}
\|\mathbf{z}_S\|_q \leq \rho \|\mathbf{z}_{\overline S}\|_1 + \tau \|\mathbf{A} \mathbf{z}\|, \quad \forall \mathbf{z} \in \mathbb{C}^N.
\end{equation}
\end{defn}
Moreover, we recall the well-known \emph{restricted isometry property} (also known as ``RIP''), where the sensing matrix is required to behave similarly to an isometry when its action is restricted to the set of sparse vectors \cite{Candes2006}.

\begin{defn}[Restricted isometry property]
The $s^{th}$ restricted isometry constant $\delta_s$ of a matrix $\mathbf{A} \in \mathbb{C}^{m \times N}$ is the smallest constant $\delta \geq 0$ such that 
\begin{equation}
(1 - \delta)\|\mathbf{z}\|_2^2 \leq \|\mathbf{A} \mathbf{z} \|_2^2 \leq (1 + \delta)\|\mathbf{z}\|_2^2, \quad \forall \mathbf{z} \in \Sigma^N_s.
\end{equation} 
The matrix $A$ has the \emph{Restricted Isometry Property} (RIP) of order $s$ if $0 < \delta_s < 1$.
\end{defn}

It is well-known that $\delta_{2s} < 4/\sqrt{41}$ is a sufficient condition for the $\ell^2$-robust NSP to hold \cite[Theorem 6.13]{Foucart2013} (we decide to use this condition and not that presented in \cite{Cai2014} for ease of exposition).

\subsection{Bounded orthonormal systems}
\label{sec:BOS}
Our analysis focuses on the case of measurement matrices $\mathbf{A}$ arising from random sampling from a \emph{Bounded Orthonormal System} (BOS) \cite{Rauhut2010}. Some significant examples of random sampling from a BOS are the subsampled Fourier transform, nonharmonic Fourier measurements, and random evaluation of orthogonal polynomials. We will discuss these case studies in more detail in Section~\ref{sec:casestudies}. 

\begin{defn}[Random sampling from a BOS]
\label{def:BOS}
Let $\mathcal{D} \subseteq \mathbb{R}^d$ be endowed with probability measure $\nu$. Then, a set $\Phi=\{\phi_1,\ldots,\phi_N\}$ of complex-valued functions on $\mathcal{D}$  is called a \emph{Bounded Orthonormal System} (BOS) with constant $K$ if, for every $j,k \in [N]$, it holds 
$\int_{\mathcal{D}} \phi_j(\tau) \overline{\phi_k(\tau)} d \nu(\tau) = \delta_{jk}$ and $\|\phi_j\|_\infty := \sup_{\tau\in\mathcal{D}}|\phi_j(\tau)| \leq K$. 
Moreover, given $m$ independent random variables $\tau_1,\ldots,\tau_m$, distributed according to $\nu$, we define the random sampling matrix $\mathbf{A}\in\mathbb{C}^{m \times N}$ associated with a the BOS $\Phi$ as
\begin{equation}
A_{kj} := m^{-\frac12}\phi_j(\tau_k), \quad \forall k\in[m],\; j\in[N].
\end{equation}
\end{defn}

A crucial property of this kind of matrices is the following. For every $\delta \in (0,1)$, assuming 
$$
m \gtrsim s \; \mathcal{L}(N,s,\delta,\varepsilon,K),
$$
the $s^{th}$ restricted isometry constant of $\mathbf{A}$ satisfies $\delta_s \leq \delta$  with probability at least $1-\varepsilon$, where the factor $\mathcal{L}(N,s,\delta,\varepsilon,K)$ is polylogarithmic in $N$ and $s$ and can be chosen in different ways  (see \cite{Rauhut2010,Foucart2013,RauhutWard,Chkifa2016}).  For example, according to \cite[Theorem 12.32]{Foucart2013}, it can be defined as follows
\begin{equation}
\label{eq:defL}
\mathcal{L} = \frac{K^2}{\delta^{2}} \max\left\{\ln^2(s)\ln\left(\frac{K^2}{\delta^{2}}s \ln(N)\right)\ln(N),\ln\left(\frac{1}{\varepsilon}\right)\right\}.
\end{equation}

\section{Robustness analysis}

In this section, we illustrate our robustness analysis. For proofs of the results stated here, we refer to \cite{Brugiapaglia2017}.

\subsection{Singular values of tall random matrices}
\label{sec:SVrandom}
Given a ``tall'' matrix $\mathbf{M} \in \mathbb{C}^{N \times m}$, we denote and sort its singular values as follows
\begin{equation}
s_1(\mathbf{M}) \geq s_2(\mathbf{M}) \geq \cdots \geq s_m(\mathbf{M}) = s_{\min}(\mathbf{M}). 
\end{equation} 
First, we provide a robust error estimate for QCBP as defined in (\ref{eq:QCBP}) assuming the $\ell^2$-robust NSP. This estimate depends on the minimum singular value of the Hermitian conjugate of the sensing matrix $\mathbf{A}$.
\begin{prop}[Robust error estimate based on the NSP]
\label{prop:l2robNSP=>BProbustness}
Let $\mathbf{A}\in\mathbb{C}^{m \times N}$ be a matrix of rank $m$ that satisfies the $\ell^2$-robust NSP of order $s$ with constants $0<\rho<1$ and $\tau >0$ with respect to $\|\cdot\|_2$ (see Definition~\ref{def:NSP}). Then,  the following error estimate holds for problem (\ref{eq:QCBP})
$$
\|\widehat{\mathbf{x}}(\eta) - \mathbf{x}\|_2 
\lesssim
\frac{\sigma_s(\mathbf{x})_1}{\sqrt{s}}  +  \eta + \sqrt{\frac{m}{s}}\frac{\max\{\|\mathbf{n}\|_2 - \eta,0\}}{\min\{s_{\min}(\sqrt{\frac{m}{N}}\mathbf{A}^*),1\}},
$$
where the hidden constant depends on $\rho$ and $\tau$. 
\end{prop}

Proposition~\ref{prop:l2robNSP=>BProbustness} shows that, in order to get robust error estimates for QCBP, we need to understand the asymptotic behavior of $s_{\min}(\sqrt{\frac{m}{N}}\mathbf{A}^*)$. In particular, our goal is to show that $s_{\min}(\sqrt{\frac{m}{N}} \mathbf{A}^*) \approx 1$. With this aim, we employ principles and ideas from the theory of random matrices with isotropic heavy-tailed columns \cite{Vershynin2012}. 

Recall that a random vector $\mathbf{z}\in\mathbb{C}^N$ is said to be \emph{isotropic} if
$\mathbb{E}[\mathbf{z}\mathbf{z}^*]=\mathbf{I}$. Then, consider random matrices of the form
\begin{equation}
\label{eq:Mmatrix}
\mathbf{M} = [\mathbf{m}_1 | \cdots | \mathbf{m}_m]\in \mathbb{C}^{N \times m},
\end{equation}
where the columns $\mathbf{m}_j$ are independent random isotropic vectors. We introduce the  \emph{cross-coherence parameter}, defined as
\begin{equation}
\label{eq:definco}
\mu :=  \frac{1}{N^2} \mathbb{E} \max_{k \in [m]} \sum_{\ell \in [m]\setminus \{k\}} |\langle\mathbf{m}_k,\mathbf{m}_\ell\rangle|^2.
\end{equation}
This parameter controls the off-diagonal part of the Gram matrix of $\mathbf{M}$. We also define the \emph{distortion parameter} as 
\begin{equation}
\label{eq:defdist}
\xi := \mathbb{E}\max_{k \in [m]} \left|\frac{\|\mathbf{m}_k\|_2^2}{N} -1\right|.
\end{equation}
It measures how far the columns of $\mathbf{M}$ are from being normalized. We notice that if $\|\mathbf{m}_k\|_2 = \sqrt{N}$ almost surely for every $k \in [m]$, then $\xi$ vanishes. Using this two parameters, we give a generalization of \cite[Theorem 5.62]{Vershynin2012}.

\begin{thm}[Singular values of heavy-tailed matrices]
\label{thm:svheavytailedcols}
Let $\mathbf{M}$ be an $N \times m$ matrix ($N \geq m$) whose columns are independent isotropic random vectors in $\mathbb{C}^N$. Then, its singular values satisfy the following asymptotic estimate 
\begin{equation}
\label{eq:svestimate}
\mathbb{E} \max_{j \in [m]} \bigg|s_j\bigg(\frac{1}{\sqrt{N}}\mathbf{M}\bigg) - 1\bigg|
 \lesssim \xi  + \sqrt{(1+\xi)\mu \ln m},
\end{equation}
where $\mu$ and $\xi$ are defined as in (\ref{eq:definco}) and (\ref{eq:defdist}), respectively.  
\end{thm}

\begin{figure}
\centering
\includegraphics[width  = 9cm]{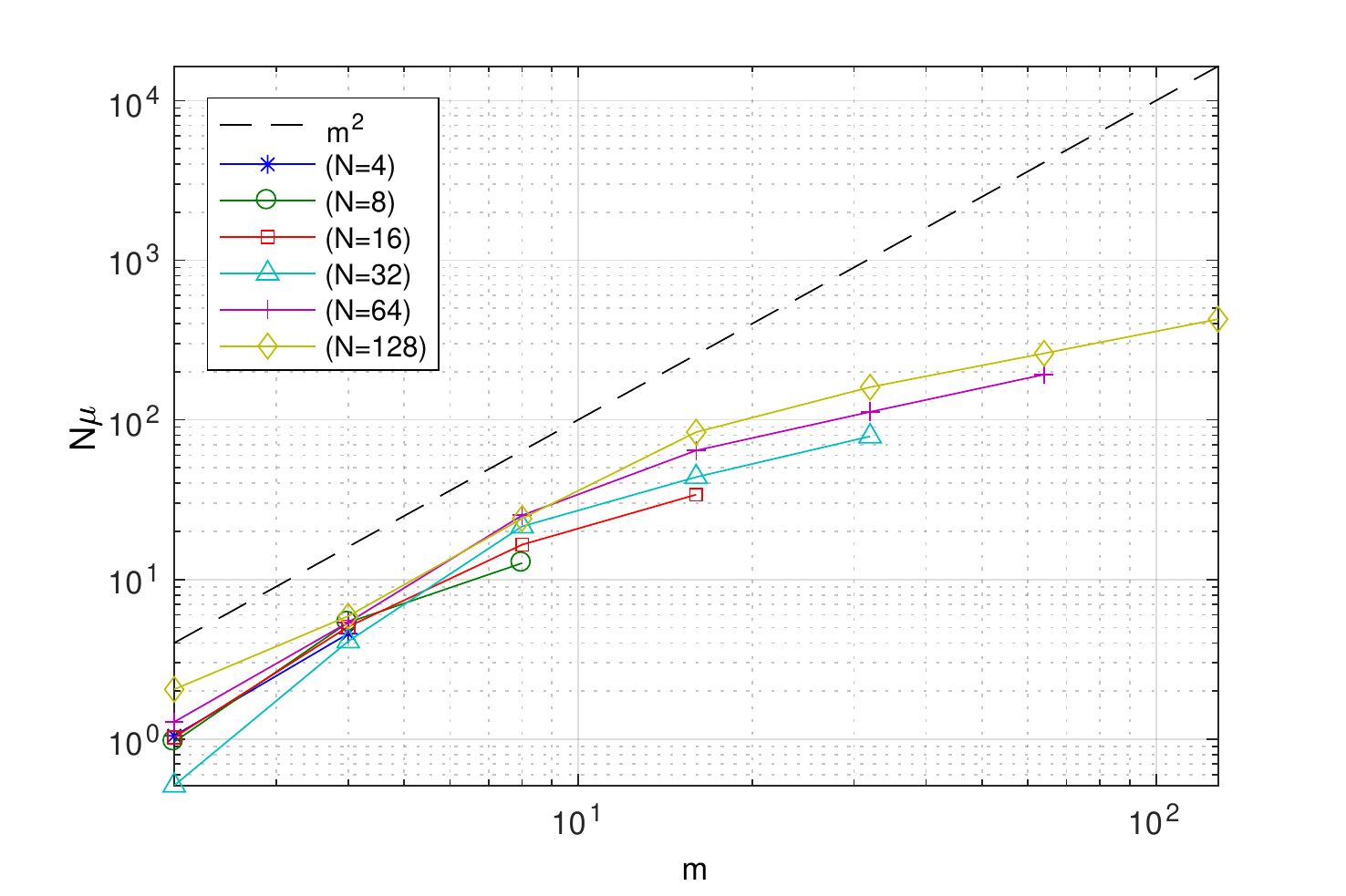}
\caption{\label{fig:incohe_param_Fourier_m} Logarithmic plot of the quantity $N\mu$, where $\mu$ is the cross-coherence parameter defined in (\ref{eq:definco}), as a function of $m$, for various values of $N$ in the case of the subsampled Fourier transform. For each value of $m$ and $N$, $\mu$ is computed by averaging over 500 random trials. The values considered are $N = 4,8,16,32,64,128$ and $m = 2^k$, with $k = 1,\ldots,\log_2(N)$. The quantity $N\mu$  is compared with the upper bound $m^2$, employed in (\ref{eq:muUB}).}
\end{figure}

We observe that a necessary condition for estimate (\ref{eq:svestimate}) to be nontrivial is that the quantity $\xi$ should be bounded uniformly in  $N$ (but not necessarily in $m$). Regarding the incoherence parameter, exploiting the isotropy of the columns of $\mathbf{M}$, and assuming $\|\mathbf{m}_j\|_2 \leq K \sqrt{N}$ for a suitable constant $K>0$ (this is always the case for $\mathbf{M} = \sqrt{m}\mathbf{A}^*$, where $\mathbf{A}$ is the random sampling matrix associated with a BOS), we can prove that
\begin{equation}
\label{eq:muUB}
N \mu \leq (Km)^2.
\end{equation}
We check the sharpness of this upper bound for moderate values of $m$ numerically in the case of the subsampled Fourier transform, where $K=1$ (see Figure~\ref{fig:incohe_param_Fourier_m}).

\subsection{Recovery error estimate}

The following theorem is the main result of the paper. It provides a robust error estimate in probability for QCBP, when  the solution $\mathbf{x}$ is very sparse and  $\mathbf{A}$ is a random sampling matrix associated with a BOS.

\begin{thm}[Robust recovery error estimate for QCBP]
\label{thm:robustrecovery}
Consider a BOS $\Phi$ with constant $K \geq 1$ and with  a distortion parameter that satisfies
\begin{equation}
\label{eq:robustrecovery:D}
\xi \leq \min\bigg\{D_1 \sqrt{\frac{m^2 \ln(m)}{N}},D_2\bigg\}, 
\end{equation}
for suitable constants $D_1,D_2$ independent of $m$ and $N$. Then, there exist constants $c,d,C,D,E>0$ and a function $\mathcal{L} = \mathcal{L}(N,s,\varepsilon,K)$ depending polylogarithmically on $N$ and $s$ such that the following holds. 
For every $N \in \mathbb{N}$ and $\varepsilon \in(0,1)$, assume that the sparsity level satisfies
\begin{equation}
\label{eq:robustrecovery:sparsity}
s \leq \frac{\varepsilon \, \sqrt{N}}{c\;\mathcal{L}(N,s,\varepsilon,K)\;\ln^{\frac12}(N)},
\end{equation}
and let $\mathbf{A}\in\mathbb{C}^{m \times N}$ be the random sampling matrix associated with $\Phi$ and with a number of measurements
\begin{equation}
\label{eq:robustrecovery:m}
m = \lceil d  \, s \, \mathcal{L}(N,s,\varepsilon,K) \rceil.
\end{equation}
Then, the following robust error estimate holds for QCBP \begin{equation}
\label{eq:robustrecovery:err}
\|\widehat{\mathbf{x}}(\eta)-\mathbf{x}\|_2 \leq \frac{C}{\sqrt{s}}\sigma_s(\mathbf{x})_1 + D \eta + E\,\mathcal{L}^{\frac12}\max\{\|\mathbf{n}\|_2 - \eta,0\},
\end{equation}
with probability at least $1-\varepsilon$. The constant $c$ depends on $D_1$, $D_2$, and $K$, whereas the constants $d$, $C$, $D$, and $E$ are universal. The function $\mathcal{L}$ can be defined as in (\ref{eq:defL}), with $\delta = 1/2$.
\end{thm}

The required relation between $s$ and $m$ is linear up to logarithmic factors, in accordance with the usual recovery error estimate in  CS. However, there are three main limitations of Theorem~\ref{thm:robustrecovery} that are worth underlining. First, the result holds for a particular sparsity regime (\ref{eq:robustrecovery:sparsity}). Essentially, fixed the failure probability $\varepsilon$,  we require $s \lesssim \sqrt{N}$ (up to logarithmic factors). Second, the error estimate (\ref{eq:robustrecovery:err}) actually depends on $N$, but this dependence is only polylogarithmic (the factor $\mathcal{L}^{\frac12}$ is due to the term $\sqrt{m/s}$ in the error estimate of Proposition \ref{prop:l2robNSP=>BProbustness}). Third, there is a linear dependence between $s$ and $\varepsilon$ in (\ref{eq:robustrecovery:sparsity}). Thus, the failure probability of the estimate is not ``overwhelmingly low''. These three issues are open problems currently under investigation.

\subsection{Applications}
\label{sec:casestudies}
To conclude, we discuss some applications of Theorem~\ref{thm:robustrecovery} to concrete examples from signal processing and high-dimensional polynomial approximation. 

\subsubsection{Fourier and Chebyshev BOSs} We discuss two examples of BOSs very popular in CS. First, we consider the subsampled Fourier transform. We have $\mathcal{D} = [N]_0$, the system is defined as
\begin{equation}
\phi_j(\tau) = \exp\big(2\pi \textrm{i} j \tau /N), \quad\forall j,k \in[N]_0,
\end{equation}
and the sampling measure $\nu$ is the uniform discrete distribution on $\mathcal{D}$. It turns out that $K=1$ and, consequently, the distortion parameter is $\xi = 0$. As a result, condition (\ref{eq:robustrecovery:D})  of Theorem~\ref{thm:robustrecovery} is satisfied.

In the case of the Chebyshev system, we consider the Chebyshev orthogonal polynomials on $\mathcal{D}=[-1,1]$, defined as
\begin{equation}
\label{eq:Cheby1}
\phi_0(\tau) \equiv 1,
\end{equation}
\begin{equation}
\label{eq:Cheby2}
\phi_{j}(\tau)=\sqrt{2}\cos(j \arccos(\tau)), \quad \forall j\in[N-1].
\end{equation} 
They form a BOS with respect to the Chebyshev measure $d \nu(\tau) = \pi^{-1}(1-\tau^2)^{-1/2}$ on $\mathcal{D}$ and with constant $K = \sqrt{2}$. By studying the normalized Christoffel function associated with this system, we estimate that the distortion parameter $\xi$ decays proportionally to $\sqrt{m/N}$. Therefore, hypothesis (\ref{eq:robustrecovery:D}) holds true and we can apply Theorem~\ref{thm:robustrecovery}.

\subsubsection{High-dimensional polynomial approximation} We assess the robustness of QCBP for polynomial approximation in high dimension \cite{Adcock2016A}. Consider the multivariate function
\begin{equation}
\label{eq:defmultidfun}
f(x) = \ln\Big(d+1+\sum_{i=1}^{d} x_i\Big), \quad x \in [-1,1]^{d}.
\end{equation}
We fix $d=10$ and we employ  the tensorized version of the Chebyshev polynomials defined in (\ref{eq:Cheby1})-(\ref{eq:Cheby2}) over $[-1,1]^{10}$ as a sparsity basis. We set maximum degree $10$ on each variable and we restrict the multi-index space to the \emph{hyperbolic cross} shape \cite{Babenko1960}. These choices leads to a total of $N=581$ degrees of freedom. We evaluate $f$ at $m = 50$  random independent sampling points $\tau_1,\ldots,\tau_m$ identically distributed according to the tensorized Chebyshev measure over $\mathcal{D}=[-1,1]^{10}$. In order to assess the robustness to unknown error, we artificially add centered gaussian noise with standard deviation $\zeta$ to the measurements, namely
\begin{equation}
\label{eq:noisymeas}
y_i = f(\tau_i) + \mathcal{N}(0,\zeta^2), \quad \forall i \in [m].
\end{equation}
In the case of nonintrusive methods for the uncertainty quantification of PDEs with random parameters, we can interpret this noise as the numerical error associated with the black-box PDE solver used to produce point-wise samples of the quantity of interest \cite{Doostan2011}.

\begin{figure}
\centering
\includegraphics[width=8cm]{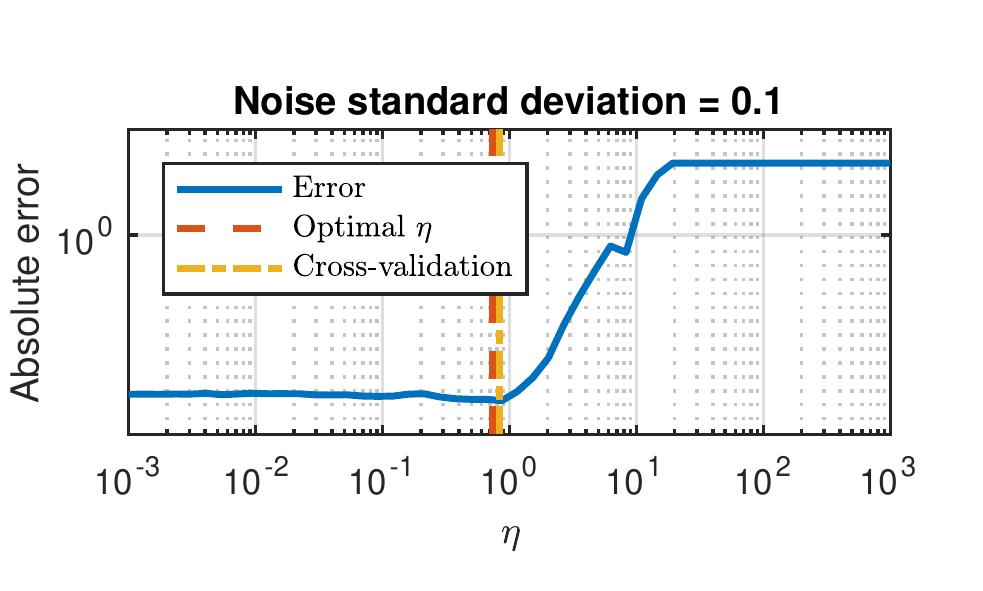}
\includegraphics[width=8cm]{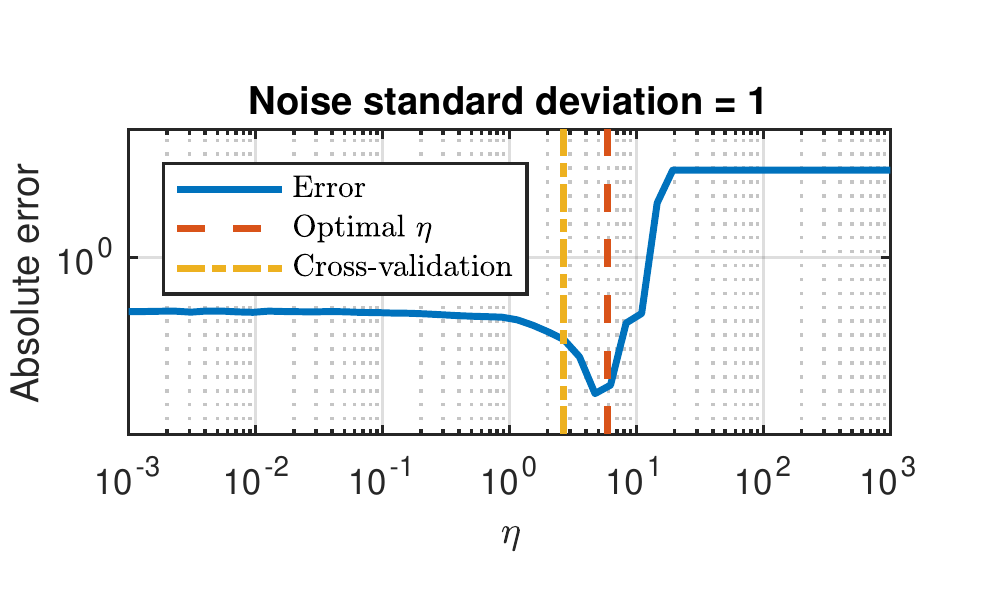}
\caption{\label{fig:multid-approx}We plot the absolute error $\|\hat{\textbf{x}}(\eta)-\mathbf{x}\|_2$ as a function of the parameter $\eta$ when QCBP is applied to the sparse polynomial approximation of the function (\ref{eq:defmultidfun}) with $d=10$. We consider 50 values of $\eta$ varying between $10^{-3}$ and $10^3$ with constant ratio and plot the results corresponding to 2 values of standard deviation $\zeta$ in (\ref{eq:noisymeas}): $\zeta = 0.1$ (top), $\zeta=1$ (bottom). We estimate the optimal value $\eta_{opt}$ as in (\ref{eq:defetaopt}). Moreover, we assess the performance of cross-validation, where $\eta$ is chosen from a grid on the interval $[10^{-2}\eta_{opt}, 10^2\eta_{opt}]$ discretized by 21 points with constant ratio. The ratio between the reconstruction samples and the validation samples is set to 3/4 (see \cite{Doostan2011} for more details). }
\end{figure}

In Figure~\ref{fig:multid-approx}, we plot the absolute error $\|\hat{\mathbf{x}}(\eta)-\mathbf{x}\|_2$ as a function of $\eta$ for different values of the standard deviation $\zeta$. The resulting curve always exhibits a global minimum. We observe that, for $\zeta =0.1$, underestimating $\eta$ is better than overestimating it. For $\zeta=1$, the minimum becomes more pronounced and underestimating $\eta$ becomes more penalizing. Recalling (\ref{eq:robustrecovery:err}) in Theorem~\ref{thm:robustrecovery}, this behavior could be justified as follows: for small values of $\zeta$, the recovery error is dominated by the term  $s^{-\frac12}\sigma_s(\mathbf{x})_1$, whereas, the more $\zeta$ gets larger, the more the term $\mathcal{L}^{\frac12}(\|\mathbf{n}\|_2 -\eta)$ becomes dominant.
   
In order to estimate the value of $\eta$ that minimizes the error, we evaluate the residual  on a reference solution $\mathbf{x}_{ref}$ computed via least-square fitting over an oversampled random grid of size $40N=23240$, i.e.,
\begin{equation}
\label{eq:defetaopt}
\eta_{opt}:=\|\mathbf{A} \mathbf{x}_{ref} - \mathbf{y}\|_2.
\end{equation}
 We compare this value with the one computed by cross-validation \cite{Doostan2011}. Both approaches are able to  approximate the minimum quite well for $\zeta=0.1$, whereas, in the case $\zeta=1$,  cross-validation slightly underperforms.

\subsubsection{Non-Cartesian Magnetic Resonance Imaging} 
\begin{figure}
\begin{center}
\begin{tabular}{cc}
\includegraphics[width = 3.2cm]{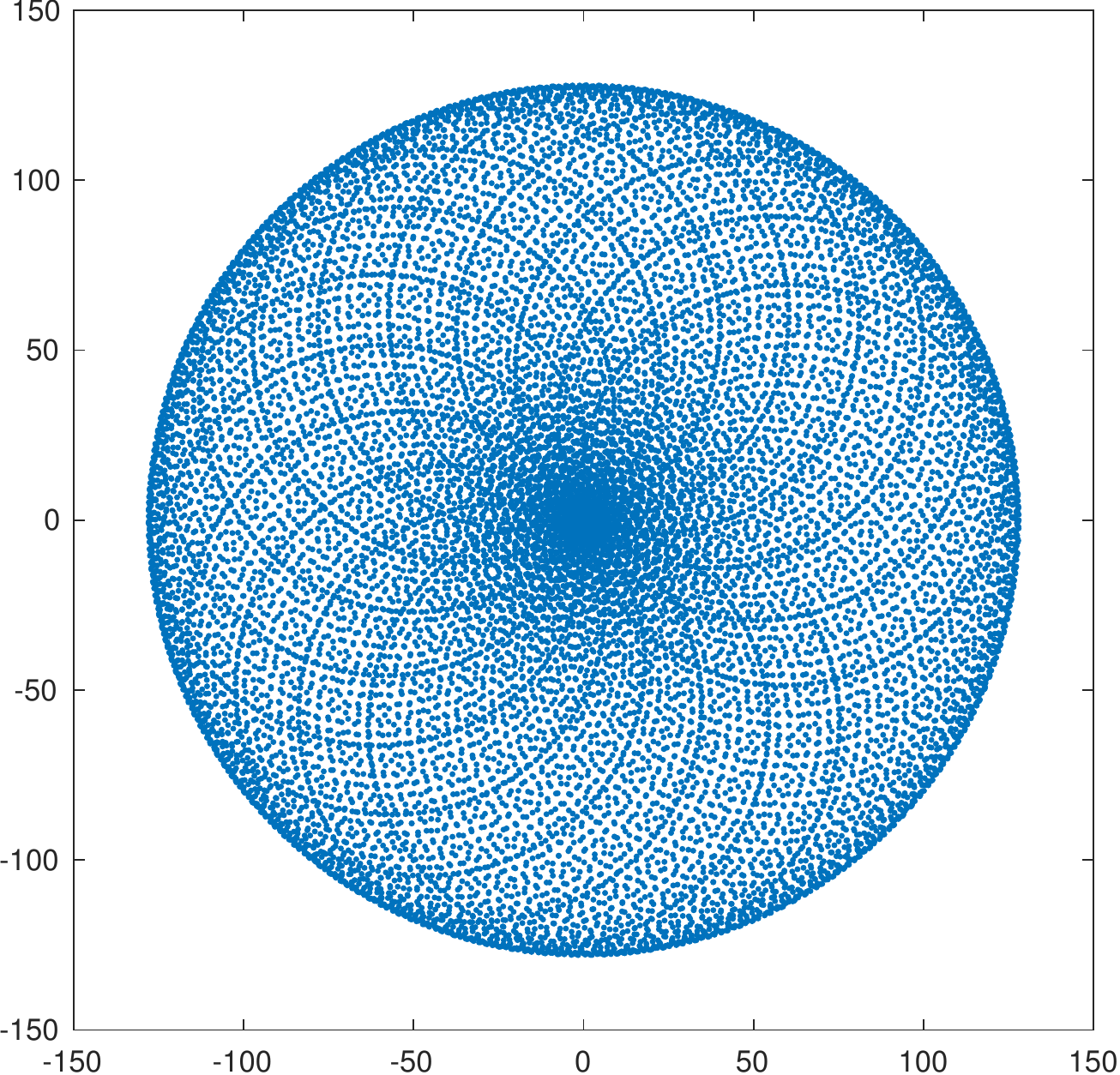} 
& 
\includegraphics[width = 3.2cm]{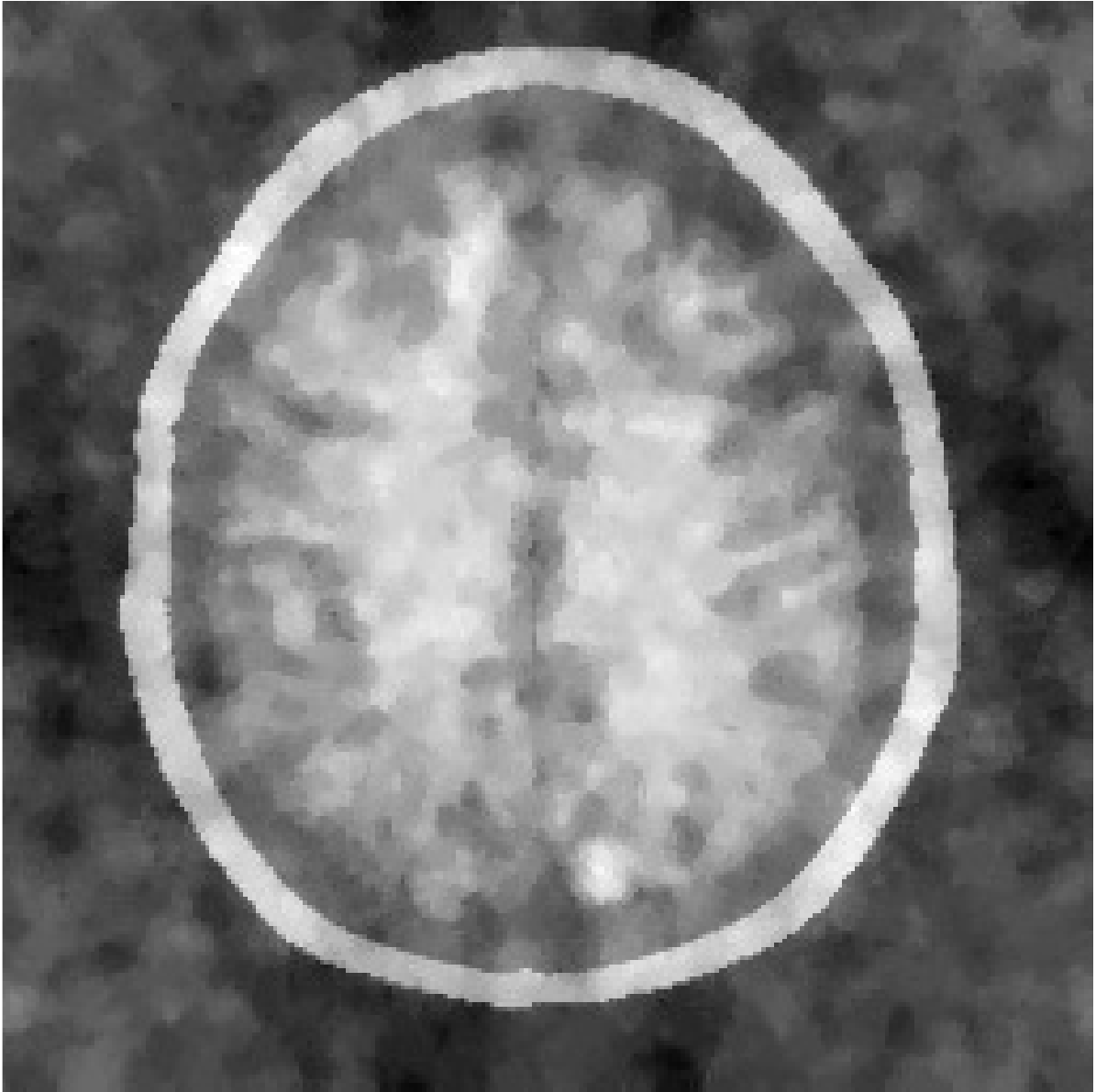} 
\\
\includegraphics[width = 3.2cm]{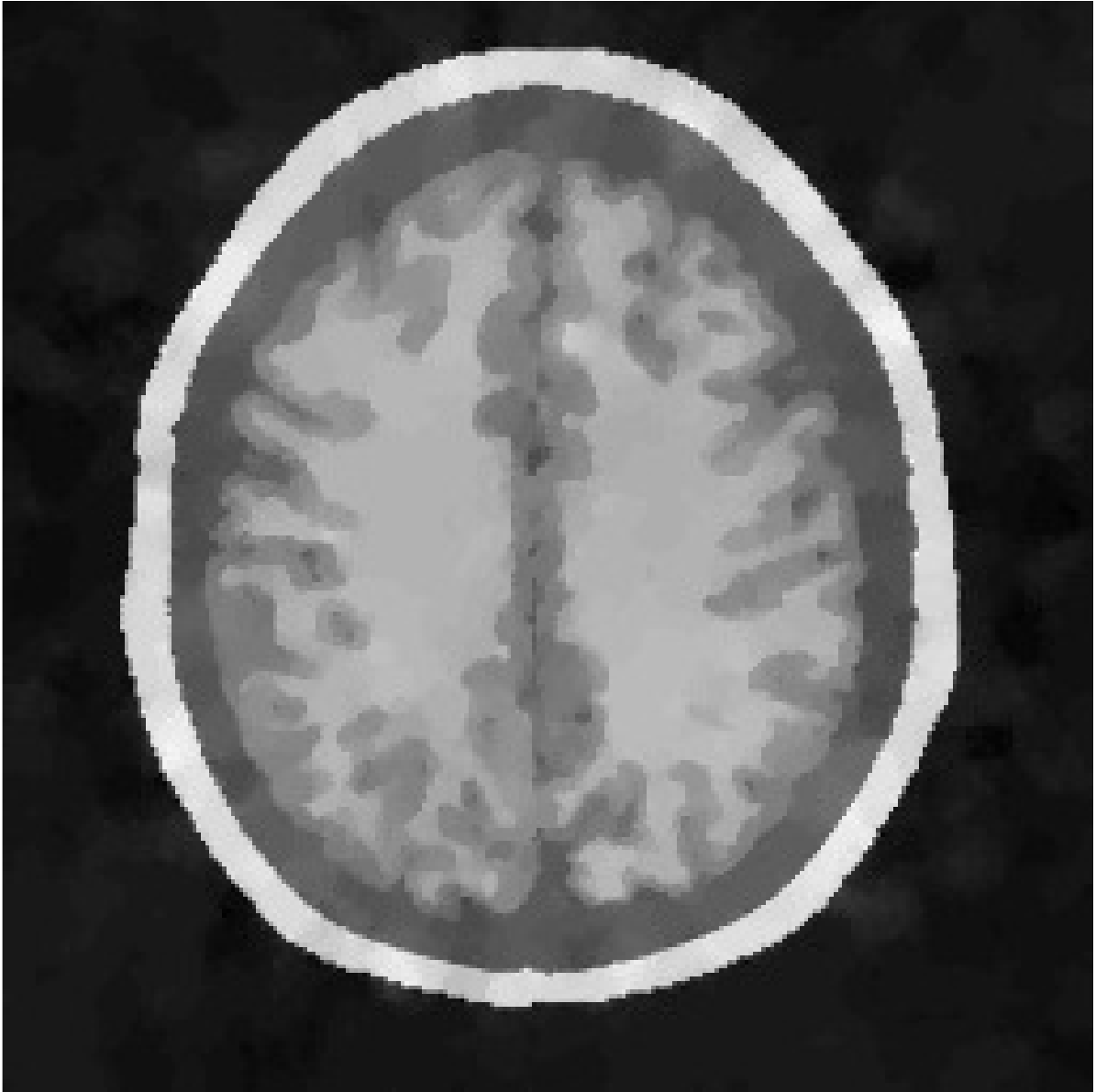} 
& 
\includegraphics[width = 3.2cm]{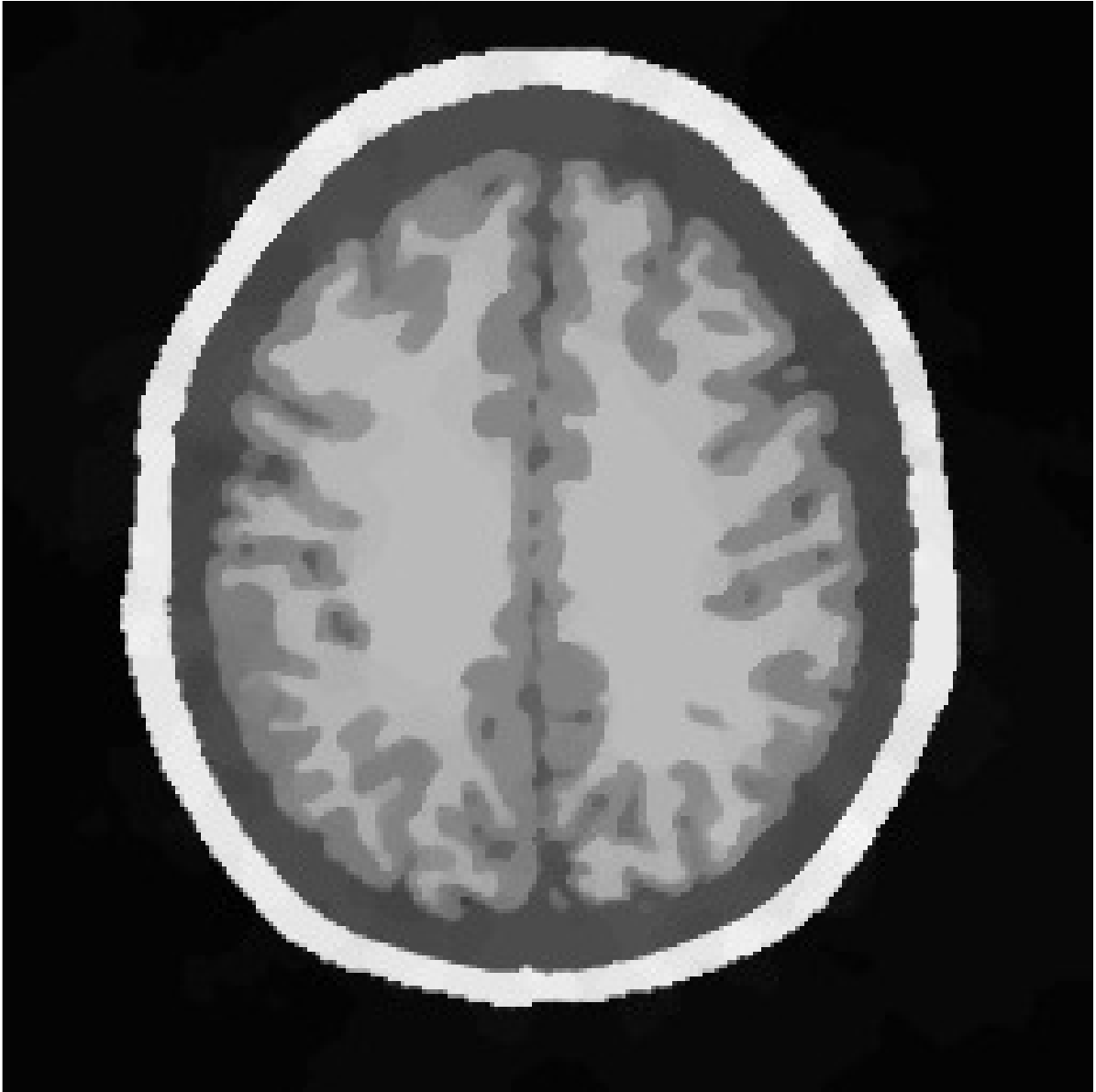}
\end{tabular}
\end{center}
\caption{\label{fig:nonCartesianMRI}The effect of gridding error on non-Cartesian MRI.  A $256 \times 256$ phantom (the McGill phantom \cite{BrainWeb}) is sampled using a rosette sampling trajectory (top left) giving $m = 15626 \approx 0.238 N$ non-Cartesian Fourier measurements.  Upfront gridding of the data is performed, and then the images are recovered using total variation minimization.  Top right: recovered image using standard nearest neighbor gridding.  The signal-to-error ratio is SER = 8.38dB.  Bottom row: recovered image using the novel fractional integer nearest neighbor gridding introduced in \cite{Adcock2016D} with parameter $n_{up} = 2$ (left) and $n_{up} = 4$ (right).  The signal-to-error ratios are SER = 14.88dB and SER = 17.84dB respectively.  Standard gridding leads to $\mathcal{O}(1)$ model error.  The fractional nearest neighbor gridding has $\mathcal{O}(1/n_{up})$ model error.}
\end{figure}

Finally, in Figure~\ref{fig:nonCartesianMRI} we give an application of this analysis to non-Cartesian MRI.  For fast reconstruction in sparse MRI, non-Cartesian data is often preprocessed by gridding it to a uniform integer grid, thus ensuring that the sampling matrix $\mathbf{A}$ can be expressed as a subsampled DFT matrix. This introduces $O(1)$ model errors, which, as seen in Figure~\ref{fig:nonCartesianMRI}, adversely affect the reconstruction.  On the other hand, the gridding strategy introduced in \cite{Adcock2016D} leads to model errors of order $1/n_{up}$, where $n_{up}$ is a user-controlled parameter.  Theorem~\ref{thm:robustrecovery} theoretically establishes the advantage of this higher-fidelity gridding, as verified in Figure~\ref{fig:nonCartesianMRI}. 


\begin{thebibliography}{1}



\bibitem{Adcock2016A}
B. Adcock, ``Infinite-dimensional {$\ell^1$} minimization and function approximation from pointwise data'' \emph{Constr. Approx.}, (to appear), 2016.

\bibitem{Adcock2016B}
B. Adcock, ``Infinite-dimensional compressed sensing and function interpolation'', \emph{Found. Comput. Math.}, (in revision), 2016.

\bibitem{Adcock2016D}
B. Adcock, R. Archibald, A. Gelb, G. Song, R. B. Platte and E. G. Walsh, ``Parameter assessment from time-dependent MR
signals using sequential imaging'', preprint, 2016.

\bibitem{Adcock2016}
B. Adcock and A.C. Hansen, ``Generalized sampling and infinite-dimensional compressed sensing,'' \emph{Found. Comput. Math.}, vol. 16, no. 5, pp. 1263--1323, 2016.

\bibitem{Babenko1960} K.I. Babenko, ``Approximation by trigonometric polynomials in a certain class of periodic functions of several variables.'' \emph{Dokl. Akad. Nauk SSSR.}, Vol. 132, no. 5, 1960.


\bibitem{BrainWeb}
BrainWeb: Simulated Brain Database. http://brainweb.bic.mni.mcgill.ca/brainweb/  

\bibitem{PhDThesisSimone}
S. Brugiapaglia, ``COmpRessed SolvING: sparse approximation of PDEs based on compressed sensing'', Ph.D.\ Thesis, Politecnico di Milano, 2016.

\bibitem{Brugiapaglia2017}
S. Brugiapaglia, B. Adcock, R.K. Archibald,
``Robustness to unknown error in sparse regularization'', \emph{in preparation}, 2017.

\bibitem{Brugiapaglia2015}
S. Brugiapaglia, S. Micheletti, and S. Perotto, ``Compressed solving: A numerical approximation technique for elliptic PDEs based on Compressed Sensing'', \emph{Comput.  Math. Appl.}, vol. 70, no. 6, pp. 1306--1335, 2015.

\bibitem{Brugiapaglia2016}
S. Brugiapaglia, F. Nobile, S. Micheletti, and S. Perotto, ``A theoretical study of COmpRessed SolvING for advection-diffusion-reaction problems'', \emph{Math. Comp.}, (to appear), 2016. 

\bibitem{Cai2014}
T. T. Cai and A. Zhang. ``Sparse representation of a polytope and recovery of sparse signals and low-rank matrices.'' \emph{IEEE Trans. Inf. Theory}, vol. 60, no. 1, pp. 122--132, 2014.


\bibitem{Candes2006}
E. J. Cand\`{e}s, J. Romberg, and T. Tao, ``Robust uncertainty principles: Exact signal reconstruction from highly incomplete frequency information,'' \emph{IEEE Trans. Inf. Theory}, vol. 52, no. 2, pp. 489--509, 2006.

\bibitem{Chkifa2016} A. Chkifa, N. Dexter, H. Tran, and C. G. Webster. ``Polynomial approximation via compressed sensing of high-dimensional functions on lower sets.'' \emph{arXiv preprint arXiv:1602.05823}, 2016.

\bibitem{Cohen2009}
A. Cohen, W. Dahmen, and R.A. DeVore, ``Compressed sensing and best k-term approximation,'' \emph{J. Amer. Math. Soc.}, vol. 22, no. 1, pp. 211--231, 2009.

\bibitem{Donoho2006}
D. L. Donoho, ``Compressed Sensing,'' \emph{IEEE Trans. Inf. Theory}, vol. 52, no. 4, pp. 1289--1306, 2006.

\bibitem{Doostan2011}
A. Doostan and H. Owhadi, ``A non-adapted sparse approximation of PDEs with stochastic inputs,'' \emph{J. Comput. Phys}, vol. 230, no. 8, pp. 3015--3034, 2011.

\bibitem{Foucart2013}
 S. Foucart and H. Rauhut, \emph{A Mathematical Introduction to Compressive Sensing}.  New York: Springer Science+Business Media, 2013.
 
\bibitem{Foucart2014}
S. Foucart, ``Stability and robustness of $\ell^1$-minimizations with Weibull matrices and redundant dictionaries,'' \emph{Linear Algebra Appl.}, vol. 441, pp. 4--21, 2014.
 
\bibitem{Guerquin-Kern2012}
M. Guerquin-Kern, L. Lejeune, K.P. Pruessmann, and M. Unser, ``Realistic analytical phantoms for parallel magnetic resonance imaging,'' \emph{IEEE Trans. Med. Imag.}, vol. 31, no. 3, pp. 626--636, 2012.
 
\bibitem{Rauhut2010}
H. Rauhut, ``Compressive sensing and structured random matrices,'' in \emph{Theoretical Foundations and Numerical Methods for Sparse Recovery}, Radon Series on Computational and Applied Mathematics, vol. 9., M. Fornasier, Ed. Berlin: de Gruyter, 2010.

\bibitem{RauhutWard}
H. Rauhut and R. Ward, ``Interpolation via weighted {$\ell_1$} minimization'', \emph{Appl. Comput. Harmon. Anal.}, vol. 40, no. 2, pp. 321--351, 2016.

\bibitem{Adcock2016C}
B. Roman, A. C. Hansen and B. Adcock, ``On asymptotic structure in compressed sensing'', arXiv:1406.4178, 2016.

\bibitem{WebsterPoly}
H. Tran, C. Webster and G. Zhang, ``Analysis of quasi-optimal polynomial approximations for parameterized {PDEs} with deterministic and stochastic coefficients'', Technical Report ORNL/TM-2015/497, Oak Ridge National Laboratory, 2015.


\bibitem{vanderBerg2008}
E. van den Berg and M.P. Friedlander, ``Probing the Pareto frontier for basis pursuit solutions,'' \emph{SIAM J. Sci. Comput.}, vol. 31, no. 2, pp. 890--912, 2008/09.



\bibitem{Vershynin2012}
R. Vershynin, ``Introduction to the non-asymptotic analysis of random matrices'' in \emph{Compressed Sensing: Theory and Applications}, Y. Eldar and G. Kutyniok, Ed. Cambridge: Cambridge University Press,  2012.

\bibitem{Wojtaszczyk2010}
P. Wojtaszczyk, ``Stability and instance optimality for Gaussian measurements in compressed sensing,'' \emph{Found. Comput. Math.}, vol. 10, no.1, pp. 1--13, 2010.



\end{thebibliography}
\end{document}